# Primes Distribution Linearity. Composite Numbers Cyclicity. Chaoticity of Primes. Twin Primes Occurrence.

*Marek Berezowski*

*Cracow University of Technology*

Corresponding address: marek.berezowski@pk.edu.pl

**Abstract.**

In this work, it was shown that all prime numbers lie on precisely defined 96 half-lines. This means that if a given number does not lie on any of the above-mentioned half-lines, then it is definitely a composite number. It was thus shown that all prime numbers satisfy a certain linear mathematical relationship. However, if a given number does not satisfy the above-mentioned relationship, then it is definitely a composite number. Additionally, an astonishing cyclical nature of numbers was demonstrated.

**Introduction.**

A certain regularity of the distribution of prime numbers was also discovered by Stanisław Ulam [Ulam, et al., 1964], [Ulam & Stein, 1967] and Martin Gardner [Gardner, 1964].

There are many articles in the scientific literature suggesting that the prime numbers are chaotic [Bershardskii, 2011], [Bogomolny, 2007], [Soundararajan, 2007], [Timberlake & Tucker, 2007], [Wolf, 2014], [Hardy & Littlewood, 1922] or has a fractal structure [Wolf, 1989]. As this article shows, both statements are true.

## 1. The Methods

We will treat successive natural numbers (*n*) as successive degrees of angle. We mark the location of these numbers in the (*x, y*) plane in a trigonometric way:

$$x(n) = n * \cos\left(n\frac{\pi}{180}\right) \qquad (1)$$

$$y(n) = n * sin\left(n\frac{\pi}{180}\right). \qquad (2)$$

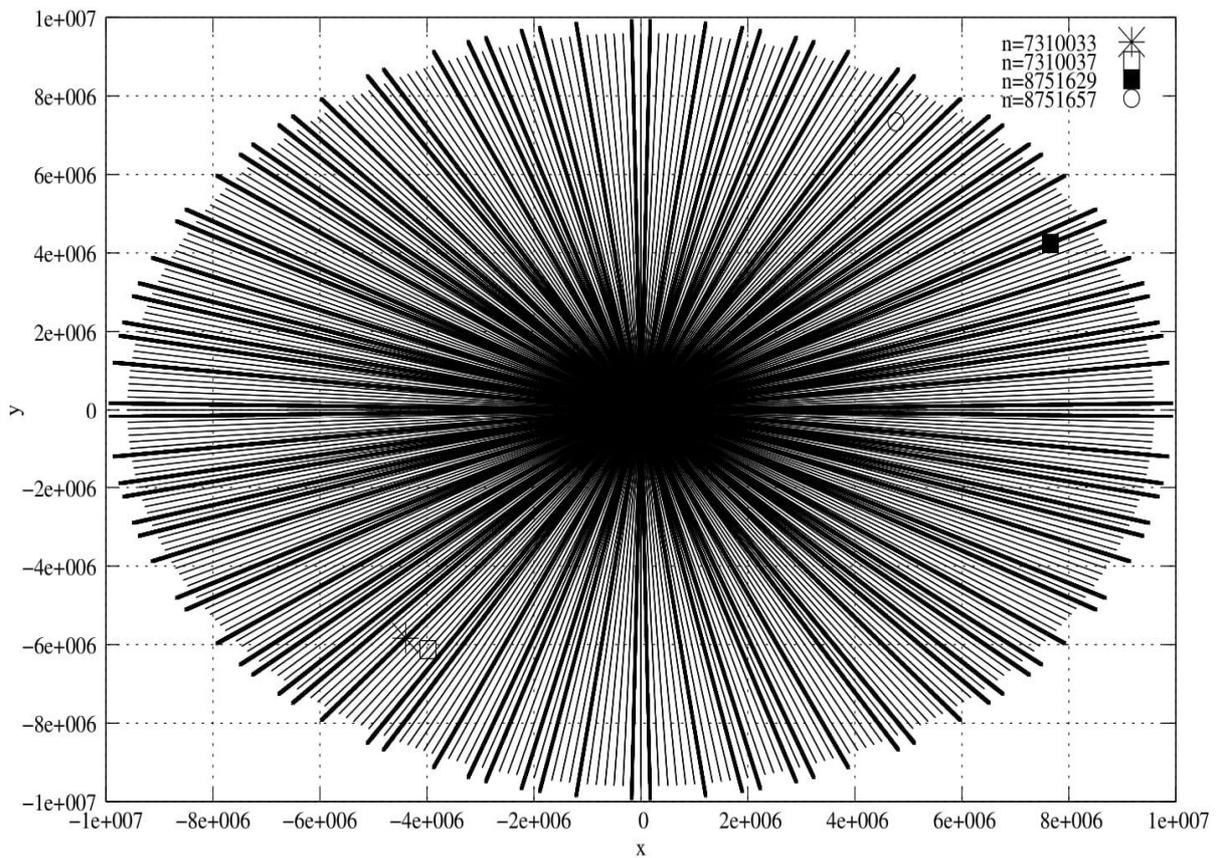

Fig. 1. Location of prime and composite numbers

On the graph, we will get a cyclically repeating set of 360 half-lines (thin and thick lines in Fig. 1). In the above graph, we will also mark the location of the prime numbers (*pn*) in the same way:

$$x(pn) = pn * \cos\left(pn\frac{\pi}{180}\right) \qquad (3)$$

$$y(pn) = pn * sin\left(pn\frac{\pi}{180}\right). \qquad (4)$$

As a result, we will obtain a cyclic bundle of 96 half-lines (thick lines on Fig. 1). 660000 prime numbers were used to draw the thick lines. Looking at Figure 1, we notice that the distribution of thick lines repeats cyclically every 30th. The first distribution contains the following sequence of 8 numbers:

$$pn_{basic} = (1, 7, 11, 13, 17, 19, 23, 29). \qquad (5)$$

For obvious reasons, this sequence does not contain primes 2, 3, and 5.

Any prime number (greater than 30) therefore satisfies the following relationship:

$$pn = pn_0 + n * 30 \qquad (6)$$

where $pn_0$ is one of the numbers in the set $pn_{basic}$, and $n$ is any natural number.

In this work, an astonishing property has been discovered. It turns out that composite numbers that do not lie on the thick rays (Fig. 1) and do not satisfy equation (6) appear cyclically in a strictly defined rhythm: 3-5-1-5-3-1-3-1, etc. This is complemented by the rhythm of numbers that lie on the thick rays and satisfy equation (6). This rhythm is: 1-2-1-2-2, etc. (Fig. 2). This means that in the case of composite numbers, three consecutive composite numbers appear first, then five, then one, then five, then three, then one, then three and finally one composite number. This cycle repeats itself indefinitely. Each of these "packages" is divided by only one number that lies on the thick rays and satisfies equation (6). For example, starting from the number 50 we have the following cycle: 50, 51, 52, *(53)*, 54, 55, 56, 57, 58, *(59)*, 60, *(61)*, 62, 63, 64, 65, 66, *(67)*, 68, 69, 70, *(71)*, 72, *(73)*, 74, 75, 76, *(77)*, 78, *(79)…* etc. The above cycle starts again from the number 80 (Fig. 2). The numbers written in brackets lie on the thick rays and satisfy equation (6), while the remaining numbers lie on the thin rays, do not satisfy equation (6) and are exclusively composite numbers. This phenomenon is also shown in Fig. 2, where the composite numbers that do not satisfy equation (6) and lie exclusively on the thick rays are marked with crosses. The circles indicate the numbers satisfying equation (6) and lying on thin rays. The individual parcels contain subsequent numbers. For example, the first parcel from the above example is the set (50, 51, 52). It should be noted that the numbers in brackets (circles in Fig. 2) are prime numbers (except numbers 77 and 91 on Fig. 2, which causes their chaotic arrangement). The remaining numbers (marked with crosses) are exclusively composite numbers. In each subsequent cycle, the situation is identical, i.e. the numbers not written in brackets are exclusively composite numbers, while the numbers in brackets can be both prime and composite numbers.

In the second case (black circles on Fig. 2), first one number satisfying (6) appears, then two, then one, then two, two, etc. Each of them is separated by at least one composite number according to the above-mentioned cycle. It should be added that in the vast majority of cases in a given cycle, the numbers in brackets are exclusively prime numbers. Using the above phenomenon, one can therefore immediately determine an arbitrarily large composite number.

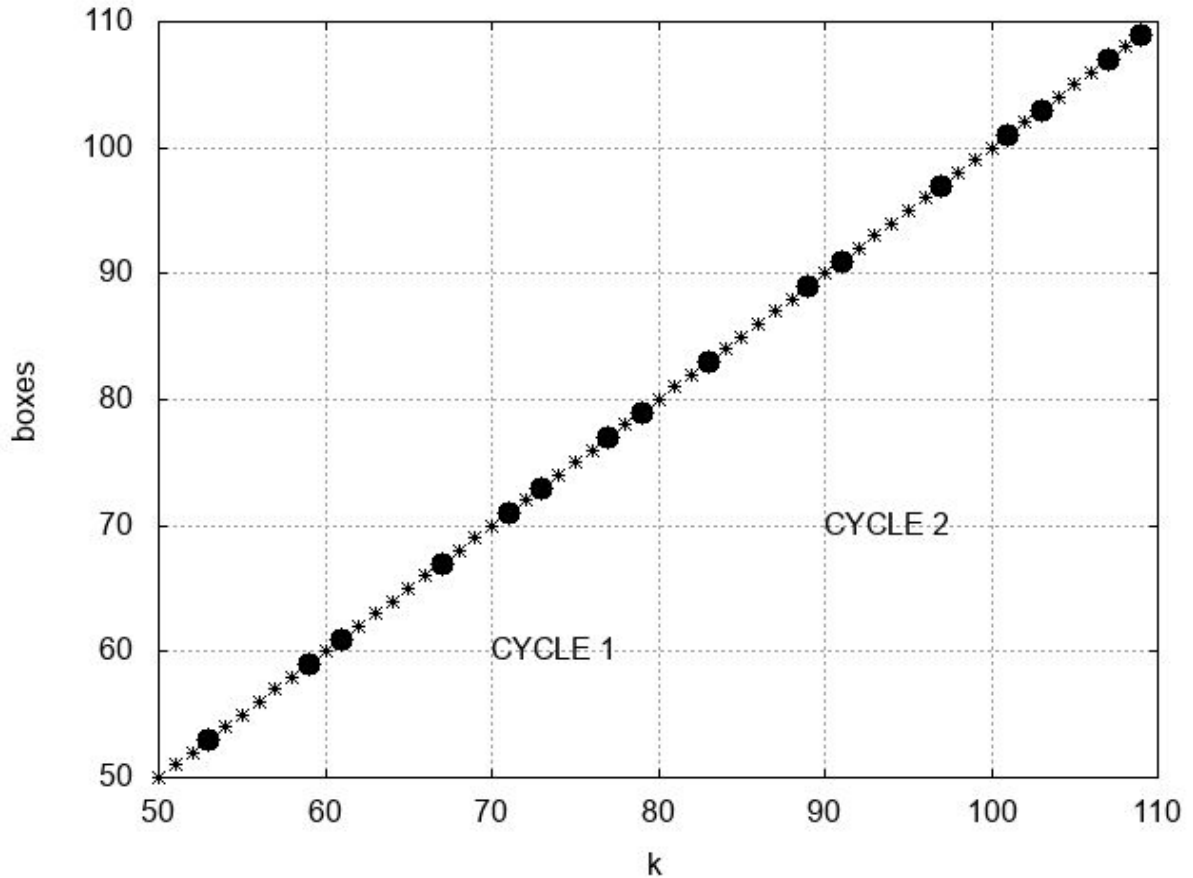

Fig. 2. Cyclicity of composite and prime numbers.

## 2. Summary.

To sum up, it should be stated that if a given natural number does not satisfy the dependence (6), it means that it is definitely a composite number. If, on the other hand, a given natural number is a prime number, it must satisfy equation (6). The graphical interpretation of this statement is that if a given natural number does not lie on any thick ray, it means that it is definitely a composite number. If, on the other hand, a given natural number is a prime number, it must lie on one of the ray lines marked in Fig. 1 with a thick line. Thus, all prime numbers are: 2, 3, 5 and all others determined by the formula (6). There are no others. It can therefore be said that the base of prime numbers is 8 numbers from the set (5). It is said that the "atoms" of natural numbers are prime numbers. The above analysis shows that the "atoms" of prime numbers are set (5). It can therefore be said that the "atoms" of all natural numbers are 8 numbers from the set (5).

An extremely interesting phenomenon is the cyclicality of the occurrence of composite numbers in the form of "boxes" according to the scheme: 3-5-1-5-3-1-3-1. Each of these boxes is separated by only one number that meets the equation (6) and lies on the thick half-line. This makes it possible to immediately determine any large composite number. Numbers that satisfy equation (6) appear in an even simpler rhythm 1-2-1-2-2. This means that in each group of 30 natural numbers there are at most 8 primes. This, in turn, means that the probability of a prime number occurring is at most 8/30. Therefore, we should look for primes only in the places marked with circles in Fig. 2. They are not present in other places. In addition, possible twin primes (*TP*) are clearly visible there (two circles separated by one star /one composite number/). In a group of 30 natural numbers there can be at most 6 of them. The probability of a twin primes occurring in the set of natural numbers is therefore at most 1/5. Starting from position 50, the values and places of their possible occurrence can be determined from the relationship: $TP=(30n+50)+k$, where $k=9,\ 21$ or $27$ and $n=0,1,2\ldots$ There are no other twin primes.

## 3. Conclusions.

Conclusion 1. All prime numbers (except 2, 3 and 5) lie on 96 thick half-lines and satisfy equation (6). This is a necessary condition for a natural number to be prime.

Conclusion 2. Natural numbers that do not lie on thick half-lines and do not satisfy dependence (6) are composite numbers. This is a sufficient condition for a given natural number to be a composite number.

## 4. Examples.

Example 1. Take the number 7310033. It is a prime number and it lies on the half-lines marked with a thick line. The coordinates of its location are (-4399287.68; -5838051.93) (asterisk).

Example 2. Take the number 7310037. It does not lie on any half-lines marked with a thick line, which means that it is definitely a composite number. The coordinates of its location are (-3981331.5; -6130712.88) (bright square). The same conclusions are reached using dependence (6), which after transformation has the following form:

$$n = \frac{pn - pn_0}{30}. \qquad (7)$$

Thus, if a given number $pn$ is prime, then there exists a number $pn_0$ in (5) for which the right-hand side of (7) is an integer. This is confirmed by the first example:

$$\frac{7310033 - 23}{30} = 243667. \qquad (8)$$

In the second example, i.e. for $pn$=7310037, there is no number $pn_0$ in (5) for which the right-hand side of (7) is an integer. So this number is a composite number.

Similar results are also obtained for all other natural numbers. Let's take some more examples and let them be the numbers 8751629 and 8751657. The first one is a prime number and lies on the fat line. Its location coordinates are (7654347.19; 4242873.93) (black square). The second number is not on the fat line. It is therefore a composite number and its location coordinates are (4766494.02; 7339757.15) (circle). In the first case, from equation (7) we get the value:

$$\frac{8751629 - 29}{30} = 291720. \qquad (9)$$

In the second case, there is no number $pn_0$ for which the right-hand side of (7) would be an integer (a composite number).

## 5. Chaos of prime numbers.

In this part the distribution of prime numbers were presented as chaotic. The appropriate spectral analysis was performed in (Berezowski, arXiv, 2023).

In Fig. 2, the numbers satisfying equation (6) are presented in the form of circles and, as mentioned above, they are arranged in the order 1-2-1-2-2 etc. This means that all prime numbers are in these circles, but it does not mean that all circles are prime numbers. So, there is no need to look for prime numbers in places not marked with circles, because they do not occur there. This definitely narrows down the size of the set in which prime numbers appear. This set is uniquely determined by equation (6). In other words, if a given number is a prime number, then it is in a circle, but the circle does not have to be a prime number. Examples in Fig. 2 are the numbers 77 and 91 (black squares). Therefore, this is only a necessary condition for the appearance of a prime number. Therefore, in order for only prime numbers to remain in the circles, some circles must be eliminated (77, 91), which means that the above order 1-2-1-2-2 will then be destroyed and the distribution of the remaining circles (and consequently, the distribution of prime numbers) becomes chaotic (Fig. 3). Black squares indicate those composite numbers that also satisfy equation (6). In the presented example, out of 16 numbers satisfying equation (6), there are only 2 such numbers: 77 and 91 (black squares). The remaining 14 numbers are prime numbers (black circles).

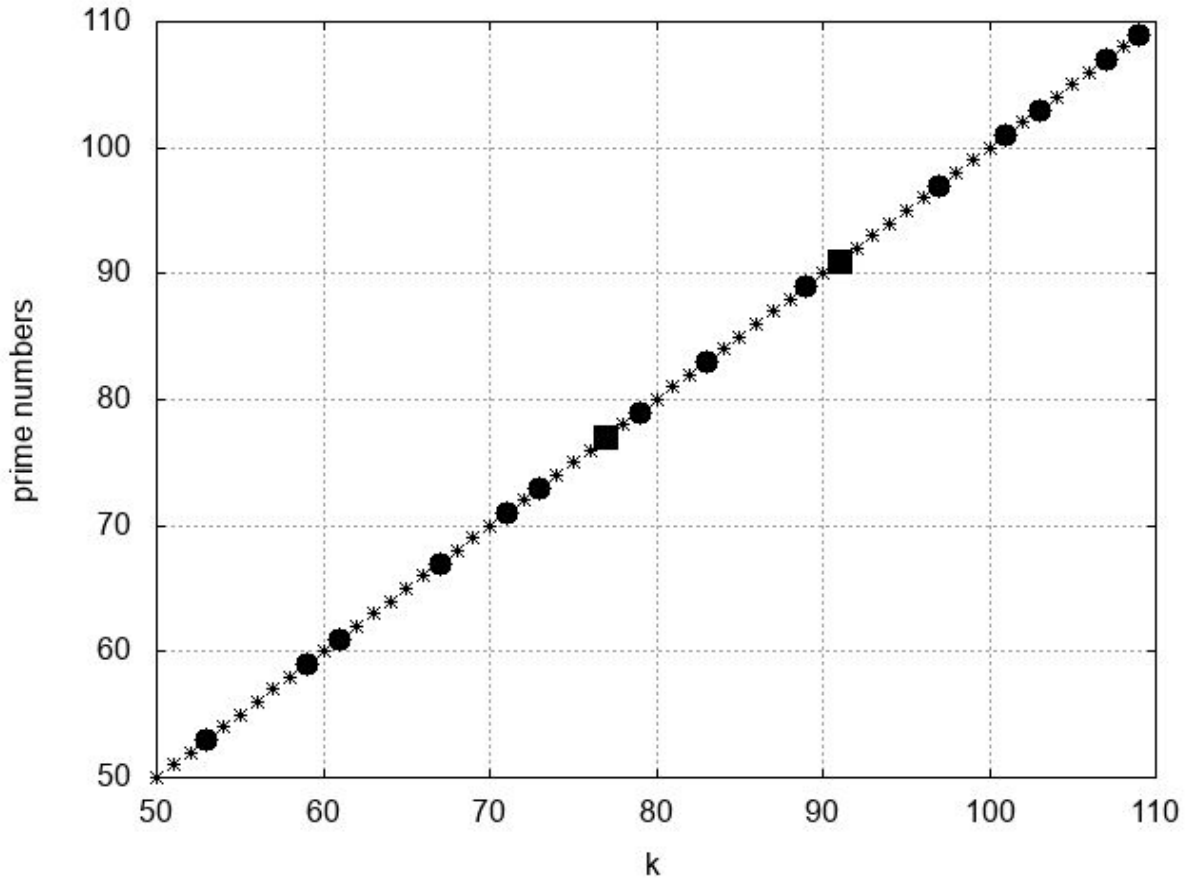

Fig. 3. Distribution of prime numbers.

## 6. Final comment.

As shown in part 1, the distribution of primes has linear character. Unlike, for example, Ulam's Spiral, all of them are located on 96 strictly defined half-lines, and all these numbers unambiguously satisfy the linear equation (6). Equally, this is based on only 8 numbers from the set (5). The presented method is much simpler and much faster, for example, than the Sieve Eratosthenes method (Havil, J., 2003). It does not require any iteration or recursion. The result is obtained in one calculation step according to equation (6) or graphically, by determining the position of the number in the diagram in Fig. 1. Let the final conclusion be that the necessary condition for a given number to be a prime number is that it lies on one of the 96 half-lines and satisfies equation (6). On the other hand, a sufficient condition for a given number to be a composite number is that it does not lie on any thick half-lines and does not satisfy equation (6). In other words, if a given natural number is to be a prime number, it must lie on one of the thick half-lines and satisfy equation (6). However, if a given number does not lie on any thick half-lines, it means that it is definitely a composite number.

So: if a given number lies on the thin half-line, then it is definitely a composite number and every prime number lies on the thick half-line. Or (which is equivalent): if a given number does not satisfy equation (6), then it is definitely a composite number and every prime number satisfies equation (6).

In the part 2 of the paper it is shown that since not all numbers satisfying equation (6) are prime numbers, the distribution of these numbers is not cyclic but chaotic. This is shown using the power spectrum in (Berezowski, arXiv, 2023). A fragment of such a distribution is shown in Fig. 3. The figure also shows exclusively the places in which twin primes can be found.